\documentclass[12pt]{article}
\usepackage{amsfonts,amsmath}
\usepackage{amssymb,latexsym}
\usepackage[dvips]{epsfig}
\usepackage{amscd}  

\title{Categorifications of the colored Jones polynomial}
\author{Mikhail Khovanov} 
\date{February 6, 2003}

\newtheorem{prop}{Proposition}
\newtheorem{theorem}{Theorem}
\newtheorem{lemma}{Lemma}

\newcommand{\oplusop}[1]{{\mathop{\oplus}\limits_{#1}}}

\begin{document}
\maketitle
\baselineskip 14pt

\def\R{\mathbb R}
\def\Z{\mathbb Z}
\def\Q{\mathbb Q}
\def\F{\mathbb F} 

\newcommand{\boldn}{\mathbf{n}}

\def\l{\lbrace}
\def\r{\rbrace}
\def\o{\otimes}
\def\lra{\longrightarrow}

\def\slt{\mathfrak{sl}(2)}

\def\mc{\mathcal} 
\def\mf{\mathfrak}
\def\sl{\mathfrak{sl}}

\def\yesnocases#1#2#3#4{\left\{
\begin{array}{ll} #1 & #2 \\ #3 & #4 
\end{array} \right. }

\newcommand{\define}{\stackrel{\mbox{\scriptsize{def}}}{=}}

\def\sbinom#1#2{\left( \hspace{-0.06in}\begin{array}{c} #1 \\ #2 \end{array}
 \hspace{-0.06in} \right)}

\def\drawing#1{\begin{center} \epsfig{file=#1} \end{center}}

\def\hsm{\hspace{0.05in}}

\def\cH{\mc{H}}
\def\cA{\mc{A}}
\def\cC{\mc{C}}
\def\cF{\mc{F}}
\def\cK{\mc{K}}
\def\cR{\mc{R}}

\begin{abstract} 
The colored Jones polynomial of links has two natural normalizations: 
one in which the $n$-colored unknot evaluates to $[n+1],$ the quantum 
dimension of the  $(n+1)$-dimensional irreducible representation of 
quantum $\slt,$ and the other in which it evaluates to $1.$ For 
each normalization we construct a bigraded cohomology theory of links
with the colored Jones polynomial as the Euler characteristic. 
\end{abstract}

\vspace{0.2in} 

\tableofcontents

\newpage


\section{The colored Jones polynomial}

\begin{figure} [htb] \drawing{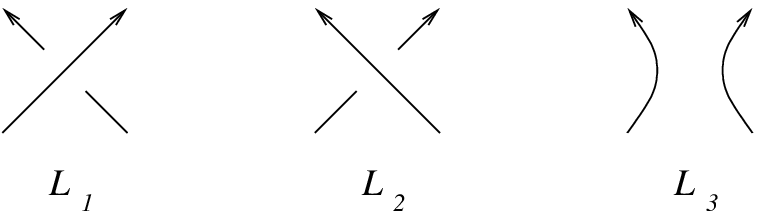}\caption{} 
  \label{skein}  
 \end{figure}

{\bf  The Jones polynomial.} 
The Jones polynomial $J(L)$ of an oriented link $L$ 
in $\R^3$ is determined by the skein relation 
 \begin{equation*} 
   q^2 J(L_1) - q^{-2} J(L_2)  = (q-q^{-1}) J(L_3)
 \end{equation*} 
for any three links $L_1,L_2,L_3$ that differ as shown
in figure~\ref{skein}, 
and its value on the unknot, which we choose to be $q+q^{-1}.$  

The Jones polynomial does not depend on the framing of link components, 
but depends slightly on their orientations. Reversing the orientation of a 
component $L'\subset L$ multiplies the polynomial by 
$q^{-6\cdot lk(L',L\setminus L')}$ 
where $lk(L',L\setminus L')$ is the linking number of $L'$ 
with its complement in $L.$ 

\vspace{0.15in} 

{\bf  The colored Jones polynomial.} We briefly recall the basics about the 
colored Jones polynomial of links (for more details consult \cite{KauLins}, 
\cite{CFS},  \cite[Sections 3,4]{KirbyMelvin}, and references therein). 
Given an oriented framed 
 link $L$ whose components are colored (marked) by non-negative 
integers (or, equivalently, by irreducible representations of 
$U_q(\mf{sl}_2),$ with integer $n$ corresponding to 
the $(n+1)$-dimensional representation $V_n$), 
the colored Jones polynomial $J_{\boldn}(L)$ takes values in 
$\Z[q,q^{-1}].$ The label $\boldn$ stands for the coloring, that is, 
for a function from the set of components of $L$ to non-negative integers. 
If all components of $L$ are colored by $1,$ the invariant is the original 
Jones polynomial of $L.$ If a component of $L$ is colored by $0,$ deleting 
this component preserves the value of the colored Jones polynomial. 

For a framed knot $K$ we denote by $J_n(K)$ the Jones polynomial of $K$ 
colored by $n.$ Thus, $J_1(K)$ is the original Jones polynomial of $K$ (and 
does not depend on the framing), while $J_0(K)=1$ for any knot $K.$ 

More generally, one could label link components by arbitrary 
finite-dimensional representations of $U_q(\slt).$ This does not give any 
extra information since the invariant is additive relative to the direct sum 
of representations, 
  \begin{equation*} 
  J_{V\oplus W}(K) = J_V(K) +  J_W(K),
  \end{equation*} 
 but allows us to express the colored Jones polynomial via the Jones 
polynomials of cables. First, the Jones polynomial 
 $J_{V\otimes W}(K)$ of a knot $K$ labelled by $V\otimes W$ 
equals the Jones polynomial of 
the 2-cable $K^2$ of $K$ with components labelled by $V$ and $W.$ 
The Jones polynomial $J_{V_1^{\otimes n}}(K)$ equals the original Jones 
polynomial of the $n$-cable $K^n$ of $K.$
 
For generic $q$ the quantum group $U_q(\slt)$ has one irreducible 
representation $V_n$ in each dimension (we consider representations 
where $q^H$ has only powers of $q$ as eigenvalues). The Grothendieck 
group of this category of representations is free abelian with 
a basis given by images $[V_n], n\ge 0 $ of all irreducible representations. 
We put square brackets around $V_n$ to distinguish a representation 
from its image in the Grothendieck group. Another basis in the 
Grothendieck group is given by all tensor powers of $V_1.$ 

Using the formula $V_n \otimes V_1 \cong V_{n+1} \oplus V_{n-1}$ and 
induction on $n,$ one checks that 
 \begin{equation} 
  [V_n] =\sum_{k=0}^{\lfloor\frac{n}{2}\rfloor} 
   (-1)^k\sbinom{n-k}{k} [V_1^{\otimes (n-2k)}]. 
 \end{equation} 
Thus, the colored Jones polynomial  $J_n(K)$ of a knot $K$ 
can be expressed via Jones polynomials of its cables: 
  \begin{equation} \label{def-knots} 
   J_n(K) =\sum_{k=0}^{\lfloor\frac{n}{2}\rfloor} 
   (-1)^k\sbinom{n-k}{k} J(K^{n-2k}). 
  \end{equation} 
This formula generalizes from knots to links, by summing  
over all connected components of a link: 
  \begin{equation} \label{def-color} 
   J_{\boldn}(L) =\sum_{\mathbf{k}=0}^{\lfloor\frac{\boldn}{2}\rfloor} 
   (-1)^{|\mathbf{k}|}\sbinom{\boldn -\mathbf{k}}{\mathbf{k}} 
   J(L^{\boldn-2\mathbf{k}}), 
   \end{equation}
where $\mathbf{k}$ is a function from the set of link components to 
non-negative integers. 

We define the colored Jones polynomial by formula (\ref{def-color}). 
Since the Jones polynomial (which appears on the right) depends 
slightly on the orientation of each cable component, we need to specify 
our choices of orientations. Each component of $L$ is oriented. 
When forming the $m$-cable of a component, if $m$ is even we orient 
half of the strands one way and the rest the opposite way. If $m$ is odd, 
we orient $\frac{m+1}{2}$ strands by way of the original orientation of the 
component, and the remaining strands the opposite way. 

\vspace{0.1in}

With this definition, $J_{\boldn}(L)$ does not depend on the orientations 
of components colored by even integers, while reversing the orientation of 
an odd-colored component $L'$ multiplies the polynomial by $q^{-6\cdot 
 lk(L',L_{\mathrm{odd}}\setminus L')}$,
where $L_{\mathrm{odd}}$ is the sublink of $L$ 
made of all odd-colored components. 

\vspace{0.1in} 

If $K$ is the 0-framed unknot,  
 \begin{equation*} 
 J_n(K)= [n+1] \define \frac{q^{n+1} - q^{-n-1}}{q-q^{-1}}, 
 \end{equation*} 
the quantum dimension of $V_n.$   

\vspace{0.2in}

\begin{figure} [htb] \drawing{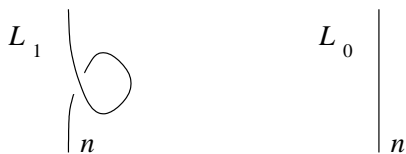}\caption{Change of framing} 
  \label{curl}  
 \end{figure}

Changing the framing of a component multiplies the colored Jones polynomial 
by a power of $q,$ as explained in the formula below and figure~\ref{curl}. 

\begin{eqnarray*} 
   J_{\boldn}(L_1) & = & q^{-2m(m+1)} J_{\boldn}(L_0)\hspace{0.15in} 
    \mathrm{if} \hspace{0.1in} n=2m, \\
   J_{\boldn}(L_1) & = & q^{-2m(m+2)} J_{\boldn}(L_0)\hspace{0.15in} 
    \mathrm{if} \hspace{0.1in} n=2m+1,
\end{eqnarray*} 
where $n$ is the color of the component. 

\vspace{0.1in}

{\bf  Reduced colored Jones polynomial.}  
Another common normalization for the colored Jones polynomial of knots is 
 \begin{equation*}  
  \widetilde{J}_n(K)\define  
  \frac{J_n(K)}{[n+1]}. 
 \end{equation*}
We will call $\widetilde{J}_n(K)$ the \emph{reduced} colored Jones 
polynomial. The reduced polynomial of the $0$-framed unknot is $1.$ We 
extend this normalization to colored links with a distinguished component, by 
  \begin{equation*}  
   \widetilde{J}_{\boldn}(L)\define \frac{J_{\boldn}(L)}{[n+1]}, 
  \end{equation*}
where $n$ is the color of the distinguished component. The reduced 
colored Jones polynomial has integral coefficients.

\section{First non-trivial examples} \label{firstex} 

{\bf  Cohomology of cables.} 
Take a framed link $L$, its cable $L^{\boldn},$ and consider the complex 
$\cC(L^{\boldn})$ and its cohomology groups $\cH(L^{\boldn})$ (we refer 
the reader to \cite{me:jones} for definitions of $\cC$ and $\cH$).   
These groups are invariants of $L,$ and depend non-trivially on its framing.  
The Euler characteristic of $\cH(L^{\boldn})$ is the Jones polynomial 
of the cable $L^{\boldn}.$ We will use the grading conventions 
of \cite{me:tangles}. 

\vspace{0.15in}

{\bf Categorification of $J_2(K).$} 
We would like to have a cohomology theory of links whose Euler characteristic 
is the colored Jones polynomial. Let's start with the first non-trivial 
example beyond the original Jones polynomial: our link is a framed knot $K$ 
colored by the three-dimensional 
representation $V_2.$ Note that $V_2$ is a direct summand of 
$V_1\otimes V_1$ and the complementary summand is the trivial representation: 
\begin{equation*}
  V_1 \otimes V_1 \cong V_2 \oplus V_0. 
\end{equation*}
This formula translates into a simple relation between 
the Jones polynomials:  
 \begin{equation} \label{2case} 
   J_2(K) = J(K^2) - 1.  
 \end{equation} 
To categorify $J(K^2)$ we take the cohomology of the 2-cable of $K.$   
Since $J_0(K)=1,$ we simply assume that its
categorification is the abelian group $\Z$ placed in bidegree $(0,0).$ 

Formula (\ref{2case}) says that to categorify $J_2(K)$ we need to 
"substract" $\Z$ from $\cH(K^2).$ This could be 
achieved by taking the cone of some map $\cH(K^2)\to \Z$ or 
of a map going in the opposite direction, from $\Z$ to  $\cH(K^2).$ 
This map should be natural. In the cohomology theory 
$\cH$ natural maps come from cobordisms between links. 
Note that $\Z$ is the cohomology of the empty link $\emptyset$ 
and there is a canonical cobordism in $\R^3\times [0,1]$ 
between $K^2$ and  $\emptyset.$ In the 2-cable $K^2$ two copies of $K$ run 
parallel next to each other, and there is a standard embedding of  an 
annulus into $\R^3$ with $K^2$ as its boundary (the two components of $K^2$ 
will be oppositely oriented if we orient our annulus and induce 
the orientation onto its boundary). Push the interior of the annulus 
into the interior of $\R^3\times [0,1]$ slighly away from the boundary 
component $\R^3\times \{ 0 \}$ that contains $K^2.$ The resulting annulus 
$S_K$ in $\R^3\times [0,1]$ is an oriented cobordism between $K^2$ and the 
empty link. Its Euler characteristic is $0$ and it induces a bidegree 
$(0,0)$ map between bigraded cohomology groups $\cH(K^2)$ and 
$\cH(\emptyset)\cong \Z.$ This map is well-defined up to overall minus sign. 
On the level of complexes, we have a map $u:\cC(K^2) \to \Z,$ well-defined 
(up to the minus sign) in the homotopy category of complexes. We define 
the complex $\cC_2(K)$ as the cone of $u,$ shifted by $[-1]$ (so that $\Z$ 
is in degree $1$), and cohomology 
groups $\cH_2(K)$ as the cohomology of $\cC_2(K).$ There is a short 
exact sequence of complexes  
 \begin{equation*} 
   0 \lra \Z[-1] \lra  \cC_2(K) \lra \cC(K^2) \lra 0 
 \end{equation*} 
giving rise to a long exact sequence of cohomology groups  
 \begin{equation*} 
    \lra  \cH_2(K) \lra \cH(K^2) \stackrel{u}{\lra} \Z \lra  
 \end{equation*} 
In this sequence at most one boundary map could be non-zero (the one from 
$\Z$ to $\cH_2^{1,0}(K)$). This map is non-zero if and only if $u$ 
is not surjective. It's easy to see, and we do it below, that $u(\cH(K^2))$ 
is a subgroup of index $1$ or $2$ in $\Z.$  

\vspace{0.2in} 

Strictly speaking, in the above discussion we should fix a plane 
diagram $D$ of $K,$ the associated diagram $D^2$ of the 2-cable $K^2,$ 
and consider the complex $\cC(D^2).$   
From two diagrams $D_0$ and $D_1$ of $K$ related by a Reidemeister 
move for framed knots, we obtain diagrams $D_0^2$ and $D_1^2$ of $K^2$ 
and a diagram of complexes and homomorphisms 
\begin{equation*} 
   \begin{CD}   
     C(D_0^2) @>{u_0}>> \Z    \\
     @VVV         @VV{\cong}V \\
     C(D_1^2) @>{u_1}>> \Z       
   \end{CD}
\end{equation*}
which commutes up to overall minus sign. Therefore, the 
cones of $u_0$ and $u_1$ are homotopy equivalent and isomorphism 
classes of groups $\cH_2^{i,j}(K)$ are invariants of $K.$ 

\vspace{0.2in}

Since $S_K$ is also a cobordism from the empty link to $K^2,$ it induces a  
 map $u': \Z\to \cC(K^2)$ and we could alternatively define 
the cohomology groups of $K$ colored by $2$ as the cohomology of the 
cone of $u'.$ 

Ideally, we would like the two resulting cohomology theories to be 
naturally isomorphic. Right off, we could see that they are if $2$ is 
invertible in the base ring. Let $k$ be a commutative ring 
where $2$ is invertible (for instance, a field of characteristic 
other than $2,$ or $\Z[\frac{1}{2}].$) Tensoring $u'$ and $u$ with $k$ we get 
maps $u': k \to \cC(K^2)\otimes_{\Z}k$ and $u:\cC(K^2)\otimes_{\Z}k \to k,$ 
and the induced maps on cohomology groups (we don't bother inventing 
different notations for these maps). The composition 
$uu': k \to k$ is the value of the invariant on the composition of the 
two cobordisms given by the annulus $S_K.$ This composition is a 
cobordism in $\R^3\times [0,1]$ between empty links, and is isotopic 
to the torus standardly embedded into $\R^4.$ Its invariant 
is $\pm 2,$ therefore, $uu'=\pm 2.$ Since $2$ is invertible, we can decompose 
 \begin{equation}\label{k-cong}  
   \cH(K^2)_k \cong \mathrm{im}(u') \oplus \mathrm{ker}(u),
 \end{equation}  
and derive isomorphisms
 \begin{equation} \label{k-cong-2} 
   \cH_2(K)_k \cong \mathrm{ker}(u) \cong \mathrm{coker}(u')
 \end{equation}  
proving that the two definitions lead to isomorphic cohomology theories
($\cH(K^2)_k$ in (\ref{k-cong}) stands for cohomology of 
$\cC(K^2)\otimes_{\Z} k$ and $\cH_2(K)_k$ in (\ref{k-cong-2}) for the 
cohomology of the total complex of 
 $0 \to \cC(K^2)\otimes_{\Z} k \stackrel{u}{\to} k \to 0 $).    

\vspace{0.1in}

The above argument implies that  $u(\cH(K^2))\subset \Z$ has index $1$ or $2.$ 

\vspace{0.1in}

\emph{Framing:} The colored Jones polynomial $J_2(K)$ depends in a 
simple way on the framing of $K.$ Change in the framing multiplies 
the polynomial by $q^{\pm 4}.$ The cohomology 
theory $\cH_2(K)$ has a more complicated dependence on the framing. Already 
when $K$ is the unknot with framing $m>0,$ the rank of  
$\cH(K^2)$ is $2 + 2m,$ and that of $\cH_2(K)$ is $1+2m.$ 
In particular, cohomology of two knots that differ by a framing are 
not just overall shifts of each other. 

In general, if $K$ is a knot and $K_1$ is obtained from $K$ by 
twisting by a large slope, 
 an overall shift superimposes $\cH_2(K)$ and $\cH_2(K_1)$ 
except for a long tail in $\cH_2(K_1)$ trailing along two 
adjacent diagonals in the bigrading plane (the tail is essentially the 
cohomology of the unknot with a large framing).   

\vspace{0.2in}

\emph{Example:} For the $0$-framed unknot $K,$ 
\begin{equation*} 
 \cH_2^{i,j}(K)\cong 
   \yesnocases{\Z}{\mbox{if } i=0 \mbox{ and }
     j\in \{ -2,0,2\}, }{0}{\mathrm{otherwise},}
\end{equation*} 
 for both definitions of $\cH_2.$ 

\vspace{0.2in}

\emph{Note:} The complex $\cC_2(K)$ and its cohomology $\cH_2(K)$ do not 
depend on the orientation of $K.$ 

\vspace{0.2in}

\begin{figure} [htb] \drawing{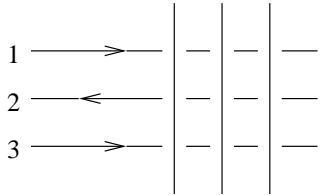}\caption{A close-up of $K^3$} 
  \label{cable3} 
 \end{figure}

{\bf Categorification of $J_3(K).$} 
Direct sum decomposition 
\begin{equation*} 
  V_1^{\otimes 3} \cong V_3 \oplus V_1 \oplus V_1 
\end{equation*} 
implies that, for a framed knot $K,$ 
\begin{equation*} 
  J_3(K) = J(K^3) - 2J(K). 
\end{equation*} 
To categorify  $J(K^3),$ we form the $3$-cable $K^3$ of $K,$ take its 
cochain complex $\cC(K^3)$ and its cohomology $\cH(K^3).$ 
Choose a plane diagram of $K$ and the corresponding diagram of the cable. 
Enumerate the components of $K^3$ by $1,2$ and $3$ so that $2$ is in the  
middle, and orient component $2$ oppositely from $1$ and $3,$ see 
figure~\ref{cable3}. 

Consider two annuli in the neighbourhood of $K$ in $\R^3,$ one with components 
$1$ and $2$ as the boundary, and the other with components $2$ and $3$ as 
the boundary. Push the interiors of the annuli into $\R^3\times [0,1]$ 
away from the boundary.  Two cobordisms from $K_3$ to $K$ result, 
inducing two maps from $\cC(K^3)$ to $\cC(K)$ and 
two maps on the cohomology, well-defined up to overall minus sign. 
We can guess that the categorification of $J_3(K)$ is the total complex  
of the bicomplex 
 \begin{equation*} 
   0 \lra \cC(K_3) \lra \cC(K) \oplus \cC(K) \lra 0. 
 \end{equation*} 
Denote the total complex by $\cC_3(K)$ and its cohomology by 
$\cH_3(K).$ We get a bigraded cohomology theory of knots with 
$J_3(K)$ as the Euler characteristic. 
Independence of $\cH_3(K)$ from the choice of planar diagram is 
straightforward. 

\vspace{0.1in}

To categorify $J_n(K)$ for arbitrary $n$ we should look for a bicomplex built 
out of complexes $\cC(K^{n-2k})$ with binomial multiplicities as in 
(\ref{def-knots}), then form the total complex and take its cohomology. This 
is done in Section~\ref{first-cat}, while the representation-theoretic 
counterpart of this construction is worked out in the next section.

%
%

\section{A resolution of an irreducible $\slt$ representation} 
\label{res-irr} 

In this section we give a homological interpretation of the 
formula 
\begin{equation} \label{resol-formula}  
 [V_n] =\sum_{k=0}^{\lfloor\frac{n}{2}\rfloor} 
  (-1)^k\sbinom{n-k}{k} [V_1^{\otimes (n-2k)}]
\end{equation} 
describing the image of $V_n$ in the Grothendieck group via those 
of tensor powers of the defining representation $V_1.$ We consider 
the $q=1$ case, so that $V_n$'s are representations of the Lie algebra 
$\slt.$ The case of generic $q$ is identical (by substituting below
"$q$-antisymmetrization" for "antisymmetrization").

The right hand side of (\ref{resol-formula}) has alternating coefficients, 
due to $(-1)^k$, and we could try to realize it as the Euler characteristic 
of a complex. 
We put the direct sum of $\sbinom{n-k}{k}$ copies of $V_1^{\otimes (n-2k)}$ 
in the $k$-th cohomological degree and look for a natural differential   
 \begin{equation*} 
   (V_1^{\otimes (n-2k)})^{\oplus \sbinom{n-k}{k}} \stackrel{d}{\lra} 
   (V_1^{\otimes (n-2(k+1))})^{\oplus \sbinom{n-k-1}{k+1}}
 \end{equation*} 
The differential reduces the number of $V_1$'s in each tensor power by 
$2.$ Since $V_1\otimes V_1 \cong V_0 \oplus V_2,$ there is a unique, up to 
scaling, surjective homomorphism $h: V_1\otimes V_1 \to  V_0,$ and it is 
simply the antisymmetrization. We will use $h$ to construct the differential.  

\begin{figure} [htb] \drawing{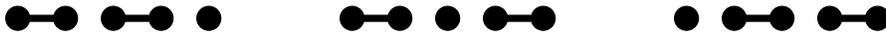}\caption{All three $2$-pairings 
 of $5$ dots} \label{pairings} 
 \end{figure}

 The binomial coefficient $\sbinom{n-k}{k}$ equals the number of 
ways to select $k$ pairs of neighbours from $n$ dots placed on a line, such 
that each dot appears in at most one pair (the 
$n=5,k=2$ example is depicted in figure~\ref{pairings}). We will call these 
$k$-pairings. 

Notice that the exponent in the tensor power of $V_1$ is the number of 
dots without a partner ("single" dots). The differential goes in the 
direction of increasing the number of pairs by $1.$ A new pair will 
consists of two adjacent   dots, and the differential should contract the two 
corresponding powers of $V_1$ into the trivial representation. 

\begin{figure} [htb] \drawing{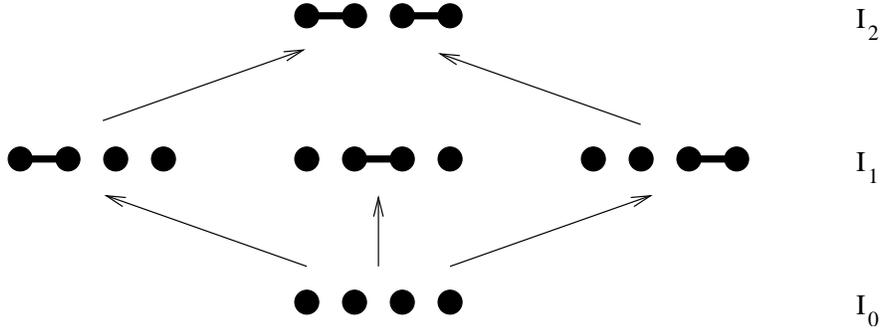}\caption{The $n=4$ case} 
 \label{graph} 
 \end{figure}

To formalize, let $I_k$ be the set of $k$-pairings of $n$ dots. 
For $s\in I_k$ let $(s)$ be the set of single dots in $s,$ and 
$V^s\define V^{\otimes (s)}$ be the tensor product 
of $V_1$'s, one for each single dot. If 
$s'\in I_{k+1}$ contains $s$ (each pair in $s$ is also a pair in 
$s'$), there is a map $h_{s',s}: V^s \to V^{s'}$ 
given by contracting the two copies of $V_1$ representing the only 
pair in $s'\setminus s.$ Consider a graph with 
with vertices--elements of $I_k,$ over all $k=0,1, \dots, 
 \lfloor \frac{n}{2} \rfloor,$ and arrows--inclusions $s\subset s'$ 
as above. For $n=4$ this graph is depicted in figure~\ref{graph}. 
Assigning $V^s$ to the vertex $s$ and the map $h_{s',s}$ to the arrow from 
$s$ to $s'$ we obtain a commutative diagram of $\slt$ representations.  

We make each square in the diagram anticommute by switching from $h_{s',s}$ 
to $(-1)^{(s,s')} h_{s',s}$ 
where $(s,s')$ is the number of pairs in $s$ to the left of the 
only pair in $s'\setminus s.$ For each $k$ take the direct sum of 
$V^s,$ for all  $s\in I_k,$ and the sum of maps as above. 
The result is a complex, denoted $C_n.$ 

\emph{Example:} The complex $C_4$ has the form 
  \begin{equation*} 
     0 \lra V_1^{\otimes 4} \lra V_1^{\otimes 2} \oplus  V_1^{\otimes 2} \oplus
     V_1^{\otimes 2} \lra V_0 \lra 0, 
  \end{equation*}
 with the differential--the sum of arrows in figure~\ref{graph}, the 
top left arrow appearing with the minus sign.  

\begin{theorem} The complex $C_n$ is acyclic in non-zero degrees 
and its degree $0$ cohomology is the irreducible $\slt$ representation 
 $V_n.$ 
\end{theorem} 

\emph{Proof:} $H^0(C_n)$ is the subrepresentation of $V_1^{\otimes n}$ 
which is the kernel of $d^0.$ This map is the sum of antisymmetrizations 
of two $V_1$'s over all pairs of neighbours. Therefore, $H^0(C_n)$ is 
isomorphic to the $n$-th symmetric power of $V_1,$ and to $V_n.$  

To prove exactness everywhere else, proceed by induction on $n$. 
The short exact sequence
 \begin{equation} \label{short-exact}  
  0 \lra C_{n-2}[-1] \lra C_n \lra C_{n-1} \otimes V_1 \lra 0 
 \end{equation} 
comes from separating all pairings of $n$ dots into two types: the 
one where the leftmost dot belongs to a pair, and the one where it does not. 
The sum of all first type selections is a subcomplex of $C_n$ isomorphic 
to $C_{n-2}$ shifted one degree to the right, the sum of all second type 
selections is a quotient complex of $C_n$ isomorphic to $C_{n-1}\otimes V_1.$ 

Induction hypothesis and the long exact sequence of (\ref{short-exact}) 
imply $H^i(C_n)=0$ for $i>1$ and exactness of 
 \begin{equation*} 
  0 \lra H^0(C_n) \lra V_1\otimes H^0(C_{n-1}) \lra H^1(C_{n-2}[-1]) 
 \lra H^1(C_n) \lra 0 
 \end{equation*} 
 Substituting $V_m$ for $H^0(C_m),$ we get a short exact sequence 
 \begin{equation*} 
   0 \lra V_n \lra V_1\otimes V_{n-1} \lra V_{n-2} \lra H^1(C_n) \lra 0 
 \end{equation*} 
 telling us that $H^1(C_n)=0.$ 

$\square$ 

\vspace{0.15in} 

\emph{Remark:}  
$C_n$ is a resolution of a simple module by a complex of semisimple modules. 
This is different from the homological algebra framework, where 
we usually resolve a module by a complex of projective, or injective, 
or flat modules. The category of finite-dimensional $\slt$ representations 
is already semisimple, all additive functors from this category are 
exact, and there is no need for resolutions from the homological 
algebra viewpoint. However, $C_n$ seems to be interesting on it own.


\section{Categorification of the colored Jones polynomial} 
\label{first-cat}

{\bf In characteristic $2.$} To avoid the sign ambiguity, in this 
subsection we work in characteristic $2.$ Let $\F_2$ be the 2-element field. 
For a diagram $D$ of a link, we denote the complex $\cC(D)\otimes_{\Z} \F_2$ 
by  $\cC(D)_2$ and its cohomology by $\cH(D)_2.$ 

\begin{figure} [htb] \drawing{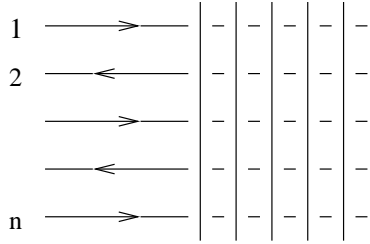}\caption{Orientations} 
 \label{orient} 
 \end{figure}

Start with a framed oriented knot $K$ and its plane diagram $D.$ It gives 
rise to diagrams $D^n$ of cables $K^n.$ The cables are oriented 
as in figure~\ref{orient}, and, in a cross-section of $D^n,$ orientations 
alternate. Choose a cross-section of $D^n$ and enumerate the strands 
from left to right from $1$ to $n$ so that component $1$ is oriented in 
the same way as $D,$ component $2$ is oppositely oriented, etc. 

For a pairing $s\in I_k$ denote by $D^s$ the cable diagram containing only 
components corresponding to single dots (not in any pair) in $s,$ and by 
$K^s$ the corresponding link (it has $n-2k$ components). 
 Given an arrow $s\to s',$ there is a canonical cobordism $S_s^{s'}$ 
from $K^s$ to $K^{s'}$ given by "contracting" the pair $s'\setminus s$ 
of neighbouring components of $D^s$ via an annulus. This cobordism 
induces a well-defined (up to homotopy) map 
$h_{s',s}: \cC(D^s)_2 \lra \cC(D^{s'})_2$ of complexes, and the 
induced map on cohomology (also denoted $h_{s',s}.$) 

Let $\cC_n(D)_2$ be the complex which in the $k$-th degree is 
the direct sum of $\cH(D^s)_2$ over all $s\in I_k.$ The differential 
is the sum of $h_{s',s}$ over all possible arrows $s\to s'.$  

Denote by $\cH_n(D)_2$ the cohomology of  $\cC_n(D)_2.$ These groups 
are bigraded, 
 \begin{equation*} 
    \cH_n(D)_2 = \oplusop{i,j\in \Z} \cH_n^{i,j}(D)_2. 
 \end{equation*} 

If diagrams $D$ and $D'$ are related by a Reidemeister move II or III, 
the complexes $\cC_n(D)_2$ and $\cC_n(D')_2$ are isomorphic. 
This implies 

\begin{prop} Isomorphism classes of bigraded groups $\cH_n(D)_2$ do not 
depend on the diagram $D$ of a framed knot $K,$ and are invariants of $K.$ 
Their Euler characteristic is the colored Jones polynomial $J_n(K).$ 
\end{prop} 

\emph{Remarks:} These groups do not depend on the orientation of $K.$ 

\vspace{0.1in} 

\emph{From knots to links.} Given a colored link $(L,\boldn),$ 
we choose its diagram $D,$ then do the above procedure for each 
component of $L.$ The result is a multi-dimensional commutative diagram 
of groups $\cH(D^{\mathbf{s}})_2$ with multiplicities--products 
of binomials, for various pairings $\mathbf{s}$ and associated 
cables $D^{\mathbf{s}}$ of $D.$ Note that $\boldn=(n_1, \dots, n_r),$ 
and the pairing $\mathbf{s}=(s_1, \dots, s_r)$ where $s_i$ is a 
pairing of $n_i$ dots (equivalently, of $n_i$ cables of the $i$-th 
component). 
We collapse this commutative 
diagram into a complex, denoted $\cC_{\boldn}(D)_2$ (in 
characteristic $2$ a commutative square is anticommutative).  
Its cohomology 
groups do not depend on $D,$ and are invariants of $L.$ 
Their Euler characteristic is $J_{\boldn}(L).$ 

\vspace{0.2in} 

\begin{figure} [htb] \drawing{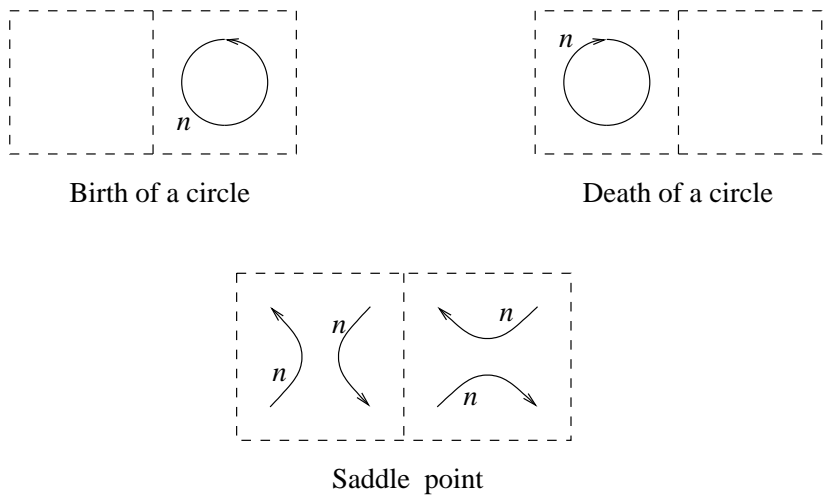}\caption{Morse moves of 
  colored link diagrams} 
 \label{morse}  
 \end{figure}

{\bf Cobordisms.} A cobordism in $\R^4$ between colored framed links 
can be represented by a sequence of its cross-sections, each  
a generic plane diagram of a colored framed link, with 
each two consequent cross-sections related by either a Reidemeister II 
or III move, or by a Morse move, see figure~\ref{morse}. 
Components involved in a saddle point move should be colored 
by the same number, see figure~\ref{morse}.   

A Reidemeister II or III move from $D$ to $D_0$ 
induces an obvious isomorphism of groups  
$\cH(D^{\mathbf{s}})_2$ and $\cH(D_0^{\mathbf{s}})_2$ 
and of complexes $\cC_{\boldn}(D)_2$ and $\cC_{\boldn}(D_0)_2.$

The unit and counit maps of $\cA^{\otimes n}$ 
induce natural maps between the complexes for  
the empty link and for the crossingless unknot diagram colored by $n.$
We assign these maps to the "birth" and "death" Morse moves. 
 
\begin{figure} [htb] \drawing{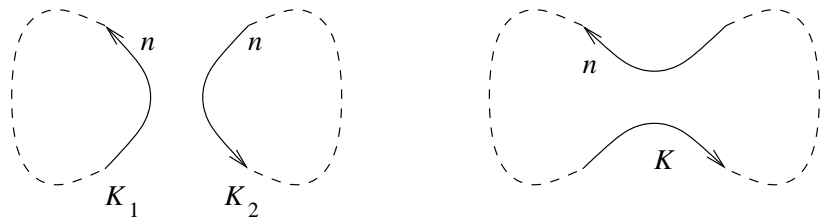}\caption{Merging two components 
 into one} 
  \label{slab}  
 \end{figure}

Suppose that diagrams $D$ and $D_0$ are related by a saddle point move. 
Consider the case when the move merges two components of $D$ (both 
labelled by $n$ and denoted $K_1$ and $K_2$) into one component $K$ 
of $D_0,$  see figure~\ref{slab}. 

We would like to come up with a natural map $\psi$ of complexes  
\begin{equation*} 
  \psi \hspace{0.05in} : \hspace{0.05in} \cC_{\boldn}(D)_2 \lra 
  \cC_{\boldn_0}(D_0)_2, 
\end{equation*} 
where $\boldn_0$ is 
the coloring of $D_0$ induced by the coloring $\boldn$ of $D.$ As an 
abelian group, $\cC_{\boldn}(D)_2$ is the direct sum of $\cH(D^{\mathbf{s}})_2$
over all pairings $\mathbf{s}$ of $\boldn.$ Given $\mathbf{s},$ 
let $s_1$ and $s_2$ be the pairings of $n$ dots which are the restrictions 
of $\mathbf{s}$ to the components $K_1$ and $K_2. $ Let $k_1$ be 
the number of pairs in $s_1$ and $k_2$ the number of pairs in $s_2.$ 

 \begin{figure} [htb] \drawing{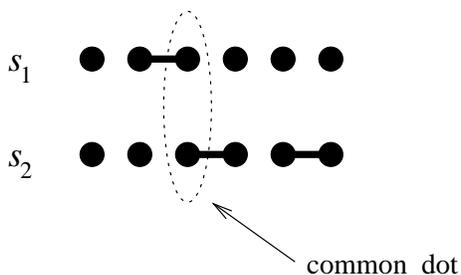}
 \caption{Two pairings with a common dot} 
  \label{common}  
 \end{figure}

\begin{figure} [htb] \drawing{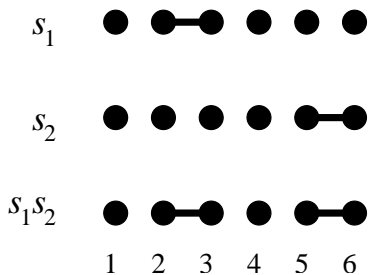}
 \caption{Disjoint union of pairings (in this example $k_1=k_2=1$)} 
  \label{disjoint}  
 \end{figure}

If pairs in 
$s_1$ and $s_2$ have at least one common dot (see figure~\ref{common}), 
we set $\psi(\cH(D^{\mathbf{s}})_2)=0.$ Otherwise, $s_1$ and $s_2$ 
are disjoint and their union $s_1s_2$ is a $(k_1+k_2)$-pairing of 
$n$ dots (see figure~\ref{disjoint}). 

To such $\mathbf{s}$ we assign a pairing 
$\mathbf{s}_0$ of the $\boldn_0$ cable of $D_0.$ This pairing 
is the same as $\mathbf{s}$ on components of $D_0$ that are unchanged 
during the saddle point move, and $s_1s_2$ on the component $K$ colored 
by $n.$ The map $\psi$ will take $\cH(D^{\mathbf{s}})_2$ to 
  $\cH(D_0^{\mathbf{s}_0})_2.$

A pair in $s_2$ connects two dots numbered $m$ and $m+1$ for some $m$ 
(in the figure~\ref{disjoint} example $m=5$). For each such pair 
we apply the operator of multiplication by $X$  
at strand $m$ of $K_1^{s_1}$ on the cohomology $\cH(D^{\mathbf{s}})_2$
(recall that $X$ is the generator of the ring $\cA,$ and $X^2=0$).  
Likewise, for a pair in $s_1$ connecting two dots numbered $t$ and $t+1$ 
(in the figure~\ref{disjoint} example $t=2$), we apply the operator 
of multiplication by $X$ at strand $t$ of $K_2^{s_2}.$ 
Denote by $\psi_1$ the product of these operators. $\psi_1$ is 
an endomorphism of $\cH(D^{\mathbf{s}})_2,$ a multiplication by 
$k_1+k_2$ copies of $X$ at certain $k_1+k_2$ strands of the cable 
$D^{\mathbf{s}}.$  

\begin{figure} [htb] \drawing{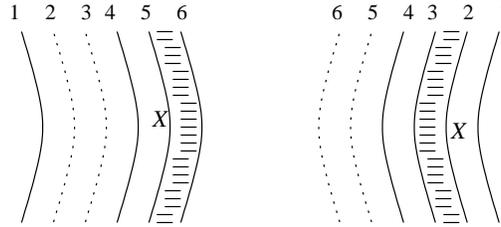}
 \caption{Multiply by $X$ at the strand $5$ of $K_1^{s_1}$ and 
  at the strand $2$ of $K_2^{s_2},$ and then contract pairs of strands 
  along annuli (these pairs are "mirrors" of dotted pairs).} 
  \label{merge2}  
 \end{figure}

For each pair in $s_2$ (connecting dots/strands numbered $m$ and $m+1$) 
consider a thin annulus whose boundary is the union of two 
strands of $K_1^{s_1}$ 
labelled $m$ and $m+1.$ Similarly, for each pair in $s_1$ 
(connecting two dots numbered $t$ and $t+1$ for some $t$) 
consider a thin annulus whose boundary is the union of two 
strands of $K_2^{s_2}$ labelled $t$ and $t+1$. Figure~\ref{merge2} 
depicts these annuli schematically for our example. The resulting 
$k_1+k_2$ annuli give rise to a cobordism from the cabled link 
with diagram $D^{\mathbf{s}}$ to a cabled link with $2(k_1+k_2)$ fewer 
components, whose diagram $D'$ can be produced by removing 
these $2(k_1+k_2)$ components (the "mirror" partners of pairs in 
$s_1$ and $s_2$) from $D^{\mathbf{s}}.$ This cobordism induces a 
map of cohomology groups from $\cH(D^{\mathbf{s}})_2$ to 
$\cH(D')_2,$ denoted $\psi_2.$  

\begin{figure} [htb] \drawing{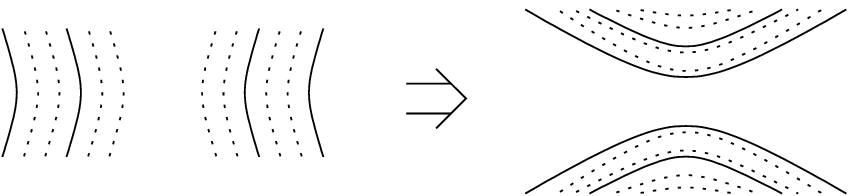}
 \caption{Cobordism which induces the map $\psi_3$ on 
  cohomology. This cobordism is a composition of saddle point 
 cobordisms between matching pairs of strands.} 
  \label{merge4}  
 \end{figure}

In the diagram $D'$ each strand in the cable $K_1^{s_1s_2}$ 
has a matching strand in the cable $K_2^{s_1s_2}.$ 
There is a canonical cobordism from $D'$ to $D_0^{\mathbf{s}_0},$ 
the composition 
of saddle point cobordisms for each pair of identically numbered 
strands in $K_1^{s_1s_2}$ and $K_2^{s_1s_2}$ (that is, one saddle 
point cobordism for each dot not in $s_1s_2$), 
for an example see figure~\ref{merge4}. Denote by $\psi_3$ 
the map of cohomology groups from $\cH(D')_0$ to $\cH(D_0^{\mathbf{s}_0})_0$ 
induced by this composition of saddle point cobordisms. 

To summarize, 
$\psi_1$ is the multiplication by a power of $X$ at certains strands 
of the cable $D^{\mathbf{s}},$ $\psi_2$ is induced by annuli contraction 
cobordisms at "mirrors" of strand pairs in $s_1$ and $s_2,$ and $\psi_3$ is 
induced by saddle point cobordisms of the remaining strands. 

Let $\psi=\psi_3 \psi_2 \psi_1:$ 
\begin{equation*} 
\psi \hspace{0.05in} : \hspace{0.1in} 
 \cH(D^{\mathbf{s}})_2 \stackrel{\psi_1}{\lra} 
  \cH(D^{\mathbf{s}})_2 \stackrel{\psi_2}{\lra} 
   \cH(D')_2 \stackrel{\psi_3}{\lra} \cH(D_0^{\mathbf{s}_0})_2. 
\end{equation*} 
Summing over all pairings 
$\mathbf{s}$ of $\boldn$ we get a homomorphism of abelian groups from 
$\cC_{\boldn}(D)_2$ to $\cC_{\boldn_0}(D_0)_2.$ 

\begin{prop} $\psi$ is a homomorphisms of complexes.  
\end{prop} 

Proof is straightforward. $\square$ 

\begin{figure} [htb] \drawing{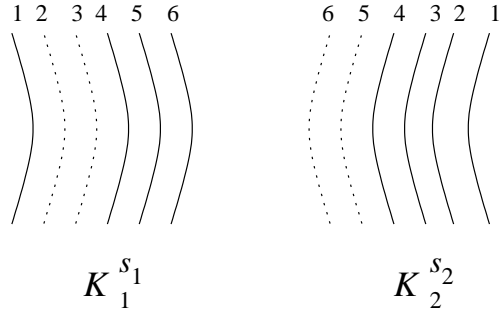}
 \caption{Part of the diagram $D^{\mathbf{s}}.$ Dotted lines show 
   strands of $D^{\boldn}$ that do not belong to 
   cables $K_1^{s_1}$ and $K_2^{s_2}.$} 
  \label{merge1}  
 \end{figure}

\begin{figure} [htb] \drawing{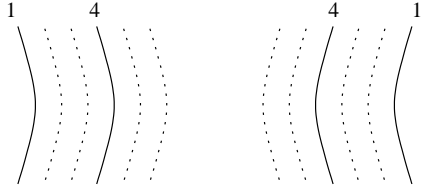}
 \caption{Diagram $D'$} 
  \label{merge3}  
 \end{figure}

Figures~\ref{merge1}-\ref{merge4} illustrate our construction 
for the case of $s_1$ and $s_2$ in figure~\ref{disjoint}.  

\vspace{0.2in}

A similar map $\psi$ can be defined for the case 
when a saddle point cobordism increases (rather than 
decreases) the number of components by $1.$ We leave its 
construction to the reader. 

\vspace{0.1in} 

Thus, to each Reidemeister and Morse move of framed colored link diagrams  
we can assign a map between the corresponding complexes 
$\cC_{\boldn}(D)_2.$ We conjecture that the induced maps on cohomology 
of these complexes give an invariant of framed colored link cobordisms. 

\vspace{0.1in} 

{\bf Categorification over a field.} Let $\F$ be a field. For a diagram $D$ of a knot $K$ 
denote by $\cC(D)_{\F}$ the complex $\cC(D)\otimes_{\Z}{\F}$ and 
its cohomology by $\cH(D)_{\F}.$  

For each arrow $s\to s'$ we have a map 
$$h_{s',s}:\cH(D^s)_{\F} \lra \cH(D^{s'})_{\F}$$ 
well-defined up to overall minus sign. For each 
square of arrows
 \begin{equation*} 
   \begin{CD}   
     s @>>> s'   \\
     @VVV         @VVV \\
     s'' @>>> s'''       
   \end{CD}
\end{equation*}
the induced square of maps $h$ either commutes or anticommutes. 

\begin{lemma}
We can always make all squares  of maps $h$ anticommute 
by changing signs in some of the maps $h_{s',s}.$
\end{lemma} 

Let us call a choice of signs in maps $h$ \emph{satisfactory} if 
all squares anticommute. 
If $h$ is satisfactory, the direct 
sum of $\cH(D^s)_{\F}$ over all possible $s$ with 
differential--the sum of $h_{s',s}$ over all possible arrows 
$s\to s'$ forms a complex. Denote this complex by $\cC_{n,h}(D)_{\F}.$ 

\begin{lemma} For any two satisfactory choices of signs 
$h', h''$ the complexes  $\cC_{n,h'}(D)_{\F}$ and  $\cC_{n,h''}(D)_{\F}$
are isomorphic. 
\end{lemma} 

Proofs of these two lemmas are left to the reader. 
From lemmas 1 and 2 we derive that the cohomology groups of
 $\cC_{n,h}(D)_{\F}$ do not depend on a satisfactory choice of signs.  
 Denote these groups by $\cH_n(D)_{\F}.$ 
Yet another exercise in sign juggling shows that 
a Reidemeister move induces an isomorphism of these cohomology groups, 
and implies 

\begin{prop} Isomorphism classes of bigraded cohomology groups $\cH_n(D)_{\F}$ 
are invariants of the framed knot $K.$ Their Euler characteristic 
is the colored Jones polynomial: 
 \begin{equation*} 
   J_n(K) = \sum_{i,j\in \Z} (-1)^i q^j \mathrm{rk}(\cH_n^{i,j}(K)_{\F}). 
 \end{equation*} 
\end{prop} 

\emph{Remark:} 
Similar considerations work over $\Z,$ with the further complication 
that we cannot pass from the complex $\cC(D^s)$ directly to its 
cohomology (since $\Z$ is not a field), and will be forced to 
work with these complexes throughout. Our diagrams 
will not be commutative or anticommutative, but rather commutative or 
anticommutative up to chain homotopy. 

\vspace{0.2in} 

{\bf Different definitions.} 
There are several competing definitions of the complex whose 
cohomology categorify the colored Jones polynomial of a link. 
Our maps between complexes assigned to cables go in the direction 
of reducing the number of strands. We could set up these maps 
to go in the opposite direction. 

From section~\ref{res-irr} we know that the algebraic 
counterpart of the complex $\cC_n(K)_{\F}$ 
has cohomology only in degree $0.$ 
Taking the hint, we could downsize our complex 
 $$ 0 \lra \cH(K^n)_{\F} \stackrel{d_0}{\lra} 
 \oplus \cH(K^{n-2})_{\F} \lra \dots $$
by restricting to the subgroup $\mathrm{ker}(d_0)$ of 
$\cH(K^n)_{\F}.$ 

The above two modifications of  $\cC_n(K)_{\F}$ can be done 
independently, and give rise to the total of four complexes (including 
 $\cC_n(K)_{\F}$). We conjecture that if $\F$ has characteristic 
$0,$ these four complexes have isomorphic cohomology groups.

\vspace{0.1in} 

\emph{Examples:} 
(a) $n=2$ case. 
We've observed in Section~\ref{firstex} that all four definitions 
give isomorphic theories over a field of characteristic different from $2.$ 

(b) $n=3$ case over $\F_2.$ The differential in  
  \begin{equation*} 
       0\lra \cH(D^3)_2 \stackrel{d_0}{\lra} 
  \cH(D)_2 \oplus \cH(D)_2 \lra 0 
  \end{equation*} 
is surjective, and  the composition $d_0 \circ d_{-1}'$ 
is a permutation, where $d_{-1}'$ is the differential in 
   \begin{equation*} 
       0\lra \cH(D)_2 \oplus \cH(D)_2 \lra \cH(D^3)_2 \lra 0. 
  \end{equation*} 
It is easy to derive that all four definitions give isomorphic theories
in this case.  

\vspace{0.1in}

In yet another approach, we could consider an 
 $S_n$-action on $\cH(K^n)$ induced by permutations (braidings) 
of strands of the cable $K^n,$ and take the 
$S_n$-invariants under this action.  
We expect that over $\Q$ this definition would give the 
same result as each the previous four. One problem with this approach 
is the projectivity of the $S_n$-action, well-defined only up to sign.


\section{Categorification of the reduced colored Jones polynomial} 

\begin{figure} [htb] \drawing{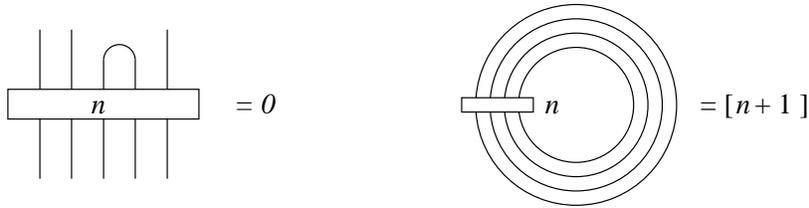}\caption{The Jones-Wenzl projector} 
   \label{joneswenzl} 
 \end{figure}

In this section we categorify the reduced colored Jones polynomial 
$\widetilde{J}_{\boldn}(L).$ Let's start with the case of a knot. 
 Given a framed colored knot $(K,n)$, the 
Jones polynomial $J_n(K)$ can be computed with the help of the Kauffman 
bracket rules and the Jones-Wenzl projector, see \cite{KauLins}. 
 The Jones-Wenzl projector is uniquely determined 
by graphical relations in figure~\ref{joneswenzl}. 

We denote the projector by $p_n$ and by $p_n'$ the projector 
divided by $[n+1].$ 

\begin{figure} [htb] \drawing{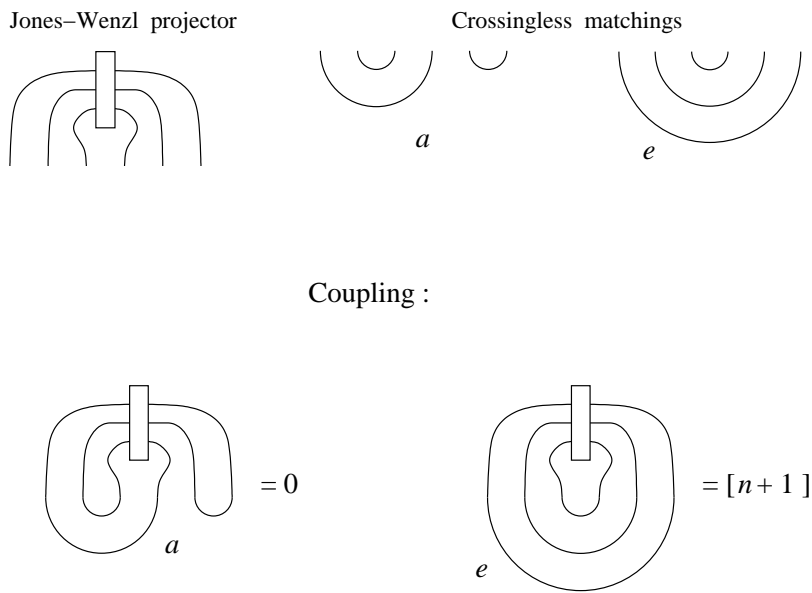}\caption{Projector, two 
crossingless matchings, and the coupling; $n=3.$} 
   \label{prmatch} 
 \end{figure}

In this section we assume familiarity with \cite{me:tangles} and, in 
particular, with the ring $H^n$ and its indecomposable left projective 
modules $P_a,$ for $a\in B^n,$ where $B^n$ is the set of crossingless 
matchings of $2n$ points. Positioning the Jones-Wenzl projector in the 
upper half of the plane
we can view it as a function from $B^n$ to $\Z[q,q^{-1}],$ by coupling 
the projector to any crossingless matching, as in figure~\ref{prmatch}.  

The Jones-Wenzl projector evaluates to zero on any crossingless matching 
except for the one denoted $e$ in figure~\ref{prmatch}, 
on which it takes value 
$[n+1].$ Therefore, $p_n'$ is a "delta-function" on $B^n$ supported on $e.$ 

\vspace{0.2in}

In \cite{me:tangles} to any crossingless matching $a\in B^n$ 
we assigned an indecomposable left projective graded $H^n$-module $P_a,$ 
and to a tangle $t$ with no bottom and $2n$ top endpoints a complex of 
projective left $H^n$-modules $\cF(t),$ well-defined in the homotopy category 
$\cK^n$ of complexes of graded $H^n$-modules. To a tangle 
$s$ with no top and $2n$ bottom endpoints we assigned a complex 
of projective right $H^n$-modules $\cF(s).$ Coupling $t$ with $s$ 
corresponds to forming the tensor product, 
$\cF(st)\cong \cF(s)\otimes_{H^n} \cF(t).$   

For each $a\in B^n$ there exists a right $H^n$-module, denoted $_a \Z,$ 
isomorphic as a graded abelian group to $\Z$ placed in degree $0,$ 
with the idempotent $1_a\in H^n$ acting as the identity on $_a \Z,$ and 
other minimal idempotents
acting by $0.$ If we were working over a field rather than over $\Z,$ this 
module would have been the simple quotient of the right projective 
module $_aP,$ while over $\Z$ the modules $_a\Z$ are substitutes for simple 
modules (for instance, they give a basis in the Grothendieck group 
of $H^n$-modules).  

We have 
 \begin{equation*}
  _a \Z \otimes_{H^n} P_b \cong \yesnocases{\Z}{\mathrm{if} \hspace{0.1in} 
  a=b,}{0}{\mathrm{otherwise.}}
 \end{equation*} 
In particular, 
 \begin{equation}
  _e \Z \otimes_{H^n} P_a \cong \yesnocases{\Z}{\mathrm{if} \hspace{0.1in}
  a=e,}{0}{\mathrm{otherwise,}}
 \end{equation}
so that $_e \Z$ has the "delta-function" behaviour when coupled to 
projective modules. Therefore, in our categorification we can interpret the 
quotient Jones-Wenzl projector $p_n'$ as the right $H^n$-module 
$_e\Z.$ The ring $H^n$ has infinite homological dimension, and $_e \Z$ 
does not have a finite length projective resolution. We do not know
any explicit construction of a (necessarily infinite) projective 
resolution of $_e \Z.$ 

\emph{Remark:} 
For a similar categorification of the Jones-Wenzl projector 
$p_n$ we should look for a graded $H^n$-module which has a filtration 
with quotient modules isomorphic to $_e \Z \{ -n\}, \hsm _e \Z \{ -n +2\}, 
\dots, \hsm _e \Z \{ n\}.$ Such a module does not exist when $n>1.$ 

We categorify the reduced colored Jones polynomial of a framed knot $K$ with 
the help of $_e \Z.$ Turn $K$ into a tangle $K_{\bullet}$ with no bottom 
and two top endpoints and form the $n$-cable 
$K^n_{\bullet}$ of this tangle. The tangle $K^n_{\bullet}$ has no bottom 
and $2n$ top endpoints. Orient it the same way we oriented cables in 
Section~\ref{first-cat}. The invariant $\cF(K^n_{\bullet})$ is a complex 
of graded projective left $H^n$-modules. Let 
  \begin{equation*} 
    \widetilde{\cC}_n(K) \define \hsm _e\Z \otimes_{H^n} \cF(K^n_{\bullet})
  \end{equation*} 
and define the reduced cohomology $\widetilde{\cH}_n(K)$ of $K$ colored by $n$ 
as the cohomology of the complex $\widetilde{\cC}_n(K).$ These cohomology 
groups are bigraded, 
  \begin{equation*} 
   \widetilde{\cH}_n(K)= \oplusop{i,j\in \Z} \widetilde{\cH}_n^{i,j}(K).
  \end{equation*} 

\begin{prop} Cohomology groups $\widetilde{\cH}_n(K)$ are invariants 
 of a framed knot $K.$ Their Euler characteristic is the reduced colored Jones 
 polynomial of $K$: 
   \begin{equation*} 
     \widetilde{J}_n(K) = \sum_{i,j\in\Z} (-1)^i q^j \mathrm{rk}
     (\widetilde{\cH}_n^{i,j}(K))
   \end{equation*} 
\end{prop} 

\vspace{0.2in} 

\emph{Remark:} 
If $K_1$ is obtained from $K_0$ by a frame change as in figure~\ref{curl}, 
we have 
 \begin{eqnarray*} 
   \widetilde{\cH}_n(K_1) & \cong & \widetilde{\cH}_n(K_0) \{-2m(m+1)\} 
  [-2m^2 ] \hspace{0.1in} \mathrm{if} \hspace{0.1in} n=2m, \\
   \widetilde{\cH}_n(K_1) & \cong & \widetilde{\cH}_n(K_0) \{-2m(m+2)\} 
  [-2m(m+1) ] \hspace{0.1in} \mathrm{if} \hspace{0.1in} n=2m+1. 
 \end{eqnarray*} 
We see that $\widetilde{\cH}_n$ depends on framing in a simpler way 
than $\cH_n.$ 

\vspace{0.2in} 

\emph{Examples:} 
\begin{enumerate} 
\item 
 When $n=0,$ reduced cohomology does not depend on the knot and 
 is isomorphic to $\Z$ placed in bidegree $(0,0).$ 
\item 
 For $n=1$ reduced cohomology was defined in \cite[Section 3]{me:patterns}.  
 $\cC(K)$ is a complex of $\cA$-modules, and tensoring it with the 
 $\cA$-module $\Z$ (where $X\in \cA$ acts trivially) gives us 
 the reduced complex $\widetilde{\cC}(K).$ 
\item 
 Reduced cohomology of the $n$-colored $0$-framed unknot is $\Z$ placed in
 bidegree $(0,0).$ 
\end{enumerate} 

For $n=1$ the relation between $\cC_n(K)$ and $\widetilde{\cC}_n(K)$
 takes the form of a short exact sequence
 \begin{equation*} 
  0 \lra \widetilde{\cC}(K)\{1\} \lra \cC(K) \lra \widetilde{\cC}(K)
 \{ -1\}  \lra 0
 \end{equation*} 
of complexes, giving rise to a long exact sequence in cohomology 
 \begin{equation*} 
   \lra \widetilde{\cH}^{i,j-1}(K) \lra \cH^{i,j}(K) \lra 
  \widetilde{\cH}^{i,j+1}(K) \lra  \widetilde{\cH}^{i+1,j-1}(K) \lra .  
 \end{equation*} 
We do not know how to relate categorifications $\cC_n(K)$ and 
$\widetilde{\cC}_n(K)$ for $n>1.$ 

\vspace{0.2in}

{\bf Categorification of the reduced colored Jones polynomial of links} 

\vspace{0.1in} 

For simplicity, in this section we switch the base ring from $\Z$ to 
the $2$-element field $\F_2.$ 

Start with a colored link $(L,\boldn)$ with a distinguished component $L'$
colored by $n.$ Turn $L$ into a tangle $L_{\bullet}$ with no bottom and 
$2$ top endpoints by cutting a segment out of $L'.$ 
Taking the $n$-cable of $L_{\bullet}$ at the component $L'$ 
gives us a tangle, denoted $L_{\bullet}^n,$ 
with no bottom and $2n$ top endpoints. 

Apply the 
construction of section~\ref{first-cat} to all components of 
$L_{\bullet}^n$ other than the $n$ components coming from $L'.$ 
The result is a complex of projective left $H^n$-modules which 
we denote $\cF_{\boldn}(L_{\bullet}).$ Let 
$$\widetilde{\cC}_{\boldn}(L)_2 \define \hspace{0.05in} _e\F_2\otimes_{H^n} 
  \cF_{\boldn}(L_{\bullet})$$
(where $_eF_2$ is the right $H_n$-module $\F_2\otimes_{\Z}(_e\Z)$),  
and define the reduced cohomology $\widetilde{\cH}_{\boldn}(L)_2$  
as the cohomology of the complex $\widetilde{\cC}_{\boldn}(L)_2.$ 

\begin{prop} Isomorphism classes of cohomology groups 
 $\widetilde{\cH}_{\boldn}(L)_2$ are invariants of a framed colored 
link $L$ with a distinguished component $L'.$ Their Euler 
characteristic is the reduced colored Jones polynomial 
$\widetilde{J}_{\boldn}(L).$  
\end{prop}


\section{Problems} 

$\quad$ 

(a) For a clean definition of cohomology groups $\cH_{\boldn}(L)$ 
it is necessary to understand the sign ambiguity of our 
cohomology theory. This ambiguity forced us into a clumsy 
definition of $\cH_{n}(K)_{\F},$ for a field $\F$ 
of characteristic other than $2,$ see section~\ref{first-cat}. 
No matter how we define  $\cH_{\boldn}(L),$ taking care of signs 
is a prerequisite for extending these homology groups to an invariant 
of (framed and colored) link cobordisms. 

\vspace{0.05in}

(b) Cohomology theory described in section~\ref{first-cat} 
depends nontrivially on the framing, and is expected to give 
rise to invariants of \emph{framed} cobordisms between framed 
links. To check the invariance, one needs a list of movie moves 
for framed cobordisms. 

\vspace{0.05in}

Once (a) and (b) have been dealt with, one can move on to the 
following problem:   

\vspace{0.05in}

(c) Give a clean definition of cohomology groups $\cH_{\boldn}(L),$ 
preferably over $\Z,$ and extend this construction to a 
functor from the category of (colored, framed, decorated) link cobordisms 
to the category of bigraded abelian groups.  

\vspace{0.1in}

Here are some other problems that we'd like to mention: 

\vspace{0.05in}

(d) Establish equivalences between various definitions of cohomology 
groups (over $\Q$) discussed at the end of section~\ref{first-cat}.  

\vspace{0.05in}

(e) Extend the categorification of the reduced colored Jones polynomial 
to an invariant of (suitably decorated) link cobordisms. 

\vspace{0.05in} 

(f) Relate categorifications of the colored Jones polynomial and 
 the reduced colored Jones polynomial. 

\vspace{0.1in} 

{\bf Acknowledgements.} I would like to thank Greg Kuperberg for 
useful discussions.



\vspace{0.1in}

Mikhail Khovanov, Department of Mathematics, University of California, 

Davis, CA 95616, mikhail@math.ucdavis.edu


\begin{thebibliography}{1}

\bibitem{CFS}
J.~S. Carter, D.~E. Flath, and M.~Saito.
\newblock {\em The classical and quantum 6$j$-symbols}.
\newblock Mathematical Notes, 43. Princeton University Press, Princeton, NJ,
  1995.

\bibitem{KauLins}
L.~H. Kauffman and S.~L. Lins.
\newblock {\em Temperley-{Lieb} recoupling theory and invariants of
  3-manifolds}.
\newblock Annals of Math. Studies, 134. Princeton University Press, 1994.

\bibitem{me:patterns}
M.~Khovanov.
\newblock Patterns in knot cohomology {I}.
\newblock arXiv math.QA/0201306.

\bibitem{me:jones}
M.~Khovanov.
\newblock A categorification of the {J}ones polynomial.
\newblock {\em Duke Math J.}, 101(3):359--426, 1999.
\newblock arXiv:math.QA/9908171.

\bibitem{me:tangles}
M.~Khovanov.
\newblock A functor-valued invariant of tangles.
\newblock {\em Algebr. Geom. Topol.}, 2:665--741, 2002.
\newblock arXiv math.QA/0103190.

\bibitem{KirbyMelvin}
R.~Kirby and P.~Melvin.
\newblock The $3$-manifold invariants of {W}itten and {R}eshetikhin-{T}uraev
  for sl(2,{C}).
\newblock {\em Invent. math.}, 105:473--545, 1991.

\end{thebibliography}
\end{document}